\documentclass{amsart}

\usepackage{amsmath,amssymb}
\usepackage{amsthm}

\newtheorem{theorem}{Theorem}
\newtheorem{proposition}{Proposition}
\newtheorem{corollary}{Corollary}
\newtheorem{lemma}{Lemma}

\theoremstyle{definition}
\newtheorem{definition}{Definition}

\theoremstyle{remark}
\newtheorem{remark}{Remark}

\newcommand {\bbC} {\mathbb{C}}

\newcommand {\bbP} {\mathbb{P}}

\newcommand {\bbR} {\mathbb{R}}

\newcommand {\bbZ} {\mathbb{Z}}

\newcommand{\id}{\mathrm{id}}

\newcommand{\bfn}{{\bf n}}

\def\L{{\mathcal L}}
\def\M{{\mathcal M}}
\def\N{{\mathcal N}}

\def\cone{{\mathop{\it Cone}}}
\def\calc{{\mathcal C}}

\def\cali{{\mathcal I}}
\def\calk{{\mathcal K}}

\def\calo{{\mathcal O}}
\def\cals{{\mathcal S}}
\def\calt{{\mathcal T}}
\def\caly{{\mathcal Y}}

\def\fuk{\mathop{{\it Fuk}}}
\def\bla{\text{--}}
\def\aut{\mathop{\rm Aut}}
\def\Pic{\mathop{\rm Pic}}
\def\End{\mathop{\rm End}}
\def\Hom{\mathop{\rm Hom}}

\def\grS{\mathop{{\rm gr}S}}
\def\torS{\mathop{{\rm tor}_\Sigma S}}
\def\qgrS{\mathop{{\rm qgr}_\Sigma S}}

\def\grR{\mathop{{\rm gr}R}}

\def\qgrR{\mathop{{\rm qgr}_\Sigma R}}

\def\grSh{\mathop{{\rm gr}\hat S}}
\def\torSh{\mathop{{\rm tor}_{\hat\Sigma}\hat S}}
\def\qgrSh{\mathop{{\rm qgr}_{\hat\Sigma}\hat S}}

\def\torPhi{\mathop{{\rm tor}_{\hat\Phi}\hat S}}
\def\qgrPhi{\mathop{{\rm qgr}_{\hat\Phi}\hat S}}
\def\Perf{{\rm Perf}}
\def\sing{{\rm sg}}
\def\Psing{\pi_{\rm sg}}

\newcommand{\mapup}[1]{\begin{picture}(30,15)(0,20)
    \put(10,30){{\small $#1$}}
    \put(2,25){\vector(1,0){26}}    
  \end{picture}
}
\newcommand{\mapback}[1]{\hspace*{-30pt}
	\begin{picture}(30,15)(0,20)
    \put(10,13){{\small $#1$}}
    \put(28,21){\vector(-1,0){26}}    
  \end{picture}
}

\newcommand{\mapdiagdown}[1]{\begin{picture}(30,30)(0,30)
    \put(10,40){{\small $#1$}}
    \put(0,45){\vector(1,-1){30}}    
  \end{picture}
}

\begin{document}

\title[Unipotence of autoequivalences]
{On the unipotence of autoequivalences of toric complete intersection
Calabi-Yau categories%
\footnote{CERN-PH-TH/2009-227}}

\author{Manfred Herbst}
\address{Institute of Mathematics, University of Augsburg, Germany}
\email{Manfred.Herbst@math.uni-augsburg.de}

\author{Johannes Walcher}
\address{PH-TH, CERN, Geneva, Switzerland}
\email{Johannes.Walcher@cern.ch}

\date{November, 2009}

\subjclass[2010]{Primary 18E30; Secondary 14J32}

\begin{abstract}
We identify a class of autoequivalences of triangulated categories of singularities
associated with Calabi-Yau complete intersections in toric varieties. Elements of this class 
satisfy relations that are directly linked to the toric data.
\end{abstract}

\maketitle
 
\setcounter{section}{-1}

\section{Introduction}

Let $X$ be a smooth projective variety over $\mathbb C$, and $D_1,\ldots, D_k$ with $k \le {\rm dim}(X) +1$ 
effective divisor classes satisfying
\begin{equation}
\label{intersect}
D_1\cap D_2\cap\cdots \cap D_k = 0.
\end{equation}
For each $i=1,\ldots, k$, we consider the autoequivalence of the bounded derived category of coherent sheaves, 
$D^b(X)$, given by tensoring with the line bundle ${\mathcal O}(D_i)$, 
i.e. for $A\in D^b(X)$,
$$
\M_i (A) = 
A\otimes \calo(D_i).
$$
Choosing for each $i$ a generic section $s_i$ of ${\mathcal O}(D_i)$, we let $\N(s_i)$ be the endofunctor
of  
$D^b(X)$ sending $A$ to
$$
\N(s_i)(A) = \cone(s_i: A \to \M_i(A)).
$$
Since as a consequence of eq.\ (\ref{intersect}), the complete intersection of $(s_1,\ldots, s_k)$ is empty, 
the associated Koszul complex $\calk(s_1,\ldots,s_k)$ is exact. Therefore,
\begin{equation}
\label{relations}
\N(s_1)\circ \N(s_2) \circ \cdots \circ\N(s_k) (A) \cong 0 
\end{equation}
for any object $A$ in $D^b(X)$. When pushed to K-theory, and denoting the image of $\M_i$ by $M_i:K(D^b(X))
\to K(D^b(X))$, we obtain the relation
$$
\prod_{i=1}^k  (M_i-{\rm id}) = 0.
$$
In this paper we obtain generalizations of these formulas in a situation coming from toric geometry and 
of interest in mirror symmetry. We now summarize the ideas involved.

Let $\bbP_{\hat\Sigma}$ be a toric Calabi--Yau variety with fan $\hat\Sigma$. Toric varieties birationally
equivalent to $\bbP_{\hat\Sigma}$ may be constructed by suitable modifcations of $\hat\Sigma$. We here 
think of modifications $\hat\Phi$ obtained by crossing a face in the secondary fan, giving rise to 
blowdowns, blowups, or flops. 
Following physics terminology, we shall refer to the various $\hat\Phi$ as different \emph{phases} of the secondary fan.
For each such modification of $\hat\Sigma$, we have an 
equivalence of triangulated categories $D^b(\bbP_{\hat\Sigma})\cong D^b(\bbP_{\hat\Phi})$. The 
corresponding functors have been described explicitly by van den Bergh \cite{vdB2002} 
and Kawamata \cite{kawamata}.
Clearly then, autoequivalences of $D^b(\bbP_{\hat\Phi})$ induce autoequivalences of 
$D^b(\bbP_{\hat\Sigma})$, and vice-versa. However, certain elements of $\aut(D^b(\bbP_{\hat\Sigma}))$ and 
relations among them are easier to see in one description than in the other, which is going to be 
the main theme of our paper. 

The primary interest, however, is in studying compact Calabi--Yau manifolds $X$ that are complete intersection 
in $\bbP_{\Sigma}$ when the latter is a smooth, toric Fano. (This is not the most general case to which
our methods apply, but one of much interest.) So let $P_1, P_2, \ldots, P_\ell$ be $\ell$ anti-effective 
divisor classes of $\bbP_{\Sigma}$, 
with $\sum_a P_a=K_{\bbP_{\hat\Sigma}}$. Choose generic sections $G_a$ 
of $\calo_{\bbP_{\hat\Sigma}}(-P_a)$ for $a=1,2,\ldots, \ell$, and consider
$$
X := \{ G_1=0, G_2=0,\ldots , G_\ell=0\} \subset \bbP_\Sigma.
$$
We recall some preliminaries about the derived category $D^b(X)$ in section \ref{Prelim}. In section 
\ref{CICY}, we give a description of $D^b(X)$ in terms of a singularity category that will be
useful for later purposes. We replace $\bbP_{\Sigma}$ with a larger combinatorial object $\bbP_{\hat\Sigma}$ 
that may be thought of as the total space of the bundle $\bigoplus_{a=1}^l \calo(P_a) \to \bbP_{\Sigma}$, 
with fibers equipped with an additional grading called \emph{R-grading}, a terminology from physics. 
Let $p_a$ be the fibre coordinate on $\calo(P_a)$. 
Then the function, a.k.a \emph{superpotential} or \emph{Landau--Ginzburg potential},
$$
W = \sum_{a=1}^l p_a G_a,  
$$
is the homological device that reduces $D^b(\bbP_{\hat\Sigma})$ 
to $D^b(X)$: The singularity
category of $W$ (namely, the triangulated category of the singularity $W:\bbP_{\hat\Sigma}\to\bbC$,
in the sense of \cite{Orlov2004,Orlov2005}) is equivalent to $D^b(X)$ (see Theorem \ref{DtoMFquotient}).

Then, for each modification $\hat\Phi$ of $\hat\Sigma$ in the above sense, one may construct
a triangulated category, $\calc_{\hat\Phi}$, as a quotient of singularity categories (Definition \ref{singcat}),
$$
\calc_{\hat\Phi} = \frac{D_{\rm sg}(\grSh)}{D_{\rm sg}(\torPhi)}.
$$
This is done in section \ref{phases}:
The category $D_{\rm sg}(\grSh)$ depends on $W$, but is independent of $\hat\Phi$, whereas the 
subcategory $D_{\rm sg}(\torPhi)$ of torsion modules depends additionally on $\hat\Phi$. 
In fact, however, all $\calc_{\hat\Phi}$ are equivalent as triangulated categories to the fixed category 
$\calc_{\hat\Sigma}\cong D^b(X)$ (Theorem \ref{allphases}). 

Working now in a fixed phase $\hat\Phi$, we associate in section \ref{main} to any toric divisor class 
$D\in \Pic(\bbP_{\hat\Sigma})$ an automorphism $\M^{D}_{\hat\Phi}$ of $\calc_{\hat\Phi}$ that can be thought of 
as tensoring with a line bundle.
We emphasize that 
the equivalences amongst the $\calc_{\hat\Phi}$ will not 
identify the $\M^{D}_{\hat\Phi}$ with each other as $\hat\Phi$ varies. (For a familiar example, see 
comments below.) 

To describe the relations analogous to (\ref{relations}) among the $\M^{D}_{\hat\Phi}$, we recall that 
the toric divisors 
$D_1,\ldots, D_n$ of $\bbP_\Sigma$ are in one-to-one correspondence with the set of one-dimensional 
cones $\Sigma(1)=(v_1,\ldots,v_n)$ of $\Sigma$. We denote the canonical sections of $\calo(D_i)$ 
by $x_i$, $i=1,\ldots n$. For $\bbP_{\hat\Sigma}$, this list is extended by (the pullbacks of) 
the divisors $P_a=:D_{n+a}$, for $a=1,\ldots,\ell$, with canonical section $p_a=:x_{n+a}$. This 
extended list is in one-to-one correspondence with the set $\hat\Sigma(1)$ of one-dimensional cones 
$(\hat v_i)_{1\le i\le n+\ell}$ of $\hat\Sigma$, and contains $\hat\Phi(1)$ for any $\hat\Phi$. 
We denote the R-grading of $x_i$ by $r_i$.
$$
r_i = \begin{cases} 0\,, & 1 \le i \le n\,, \\ 2\,, & n+1 \le i \le n+\ell\,. \end{cases}
$$

Our main result will follow straightforwardly from these definitions in section \ref{main}.
\begin{theorem}
\label{bigtheorem}
For each toric divisor $D_i$, with canonical section $x_i$ of $\calo(D_i)$,
define an endofunctor $\N_{\hat\Phi}$ of $\calc_{\hat\Phi}$ by 
$$
\N_{\hat\Phi}(x_i)(\bla) = \cone ( x_i : \bla \to \bla \otimes \calo(D_i)[r_i]).
$$
Then for each subset $\cali\subset \{1,\ldots, n+\ell\}$ such that the corresponding set of 
edges $(\hat v_i)_{i\in I}$ is not contained in any cone of $\hat\Phi$, we have the relation
$$
\bigcirc_{i\in \cali} \N_{\hat\Phi}(x_i) \cong 0.
$$
\end{theorem}
A simple consequence
is the following 
\begin{corollary}
\label{Ktheory}
Let $M^{D}_{\hat\Phi}$ be the automorphism induced by $\M^{D}_{\hat\Phi}$ on the K-theory 
$K(\calc_{\hat\Phi})$. Then for each $\cali$ as above,
$$
\prod_{i\in \cali} \bigl(M^{D_i}_{\hat\Phi}-\id\bigr) = 0
$$
\end{corollary}

Our original motivation for this work was to understand the generalization of an old observation
of Kontsevich. We let $\caly$ be the family of Calabi--Yau manifolds that is mirror to $X$ according
to Batyrev's construction. Let $B$ be the base of the family after removing the singular fibers (and
possibly more). Monodromies around loops in $B$ induce symplectic transformations (generalized
Dehn twists) that can be lifted to autoequivalences of the symplectic category (the Fukaya category
$\fuk(Y)$, where $Y$ is a generic fiber of $\caly$). Via the homological mirror symmetry (HMS) conjecture, 
$\fuk(Y)\cong D^b(X)$, one is led to expect the existence of a monodromy representation
$$
\rho: \pi_1(B) \to \aut(D^b(X)) 
$$
that has attracted some attention over the years. 

When $X$ is the quintic threefold, and $\caly$ its mirror family, 
we may model $B$ as $\bbP^1\setminus
\{0,1,\infty\}$. Note that the point at $\infty$ does not correspond to a singular threefold, but to one
with an additional automorphism of order $5$. Somewhat surprisingly, this symmetry is only 
realized projectively at the categorical level. Indeed, if $\gamma_\infty$ is a path around $\infty$, 
the categorical monodromy $\M_\infty = \rho(\gamma_\infty)\in \aut(D^b(X))$ (modulo HMS conjecture) 
satisfies the relation
\begin{equation}
\label{kontsevich}
(\M_\infty)^5 \cong (\bla)[2]
\end{equation}
To fully appreciate this remarkable relation, we rewrite this in more familiar terms using 
$\gamma_\infty^{-1} = 
\gamma_0\circ \gamma_1$. According to HMS for the quintic threefold, the monodromy $\gamma_0$ around 
$0\in \bbP^1$ corresponds to the autoequivalence $\M_0$ of tensoring with the line bundle $\calo_X(1)$, 
whilst monodromy $\gamma_1$ 
(around the conifold) is realized as twist $\calt_{\calo_X}$ by the structure sheaf \cite{ST2000}. These 
transformations of 
$D^b(X)$ can be checked to satisfy \cite[Chapter 7.1.4]{Aspinwall2004} 
\begin{equation*}
\bigl( \M_0 \circ \calt_{\calo_X}\bigr)^5 \cong (\bla)[2],
\end{equation*}
which is the way in which (\ref{kontsevich}) is often quoted.
(Some hypersurfaces in weighted projective space are treated in \cite{CK2007}.)

Our results give a uniform treatment of such relations for the general class of Calabi--Yau complete 
intersections in toric varieties. We elaborate on this point of view, together with some
other applications, in section \ref{applications}.

\subsection*{Acknowledgements}
We thank David Favero, Robert Karp, Ludmil Katzarkov, Maximilian Kreuzer, Marc Nieper-Wi\ss kirchen, and 
Dmytro Shklyarov for valuable discussions and correspondence. M.~H. was supported by the ERC Starting 
Independent Researcher Grant StG No. 204757-TQFT.

\section{Preliminaries} 
\label{Prelim}

Consider $\bbP_\Sigma$, a complete, smooth toric variety defined by a fan $\Sigma$ in the lattice $N 
\cong \bbZ^{n-k}$. Let $v_i$ for $i=1,\ldots,n$ be the generators of the one-dimensional cones of 
$\Sigma$. We briefly recall the construction of $\bbP_\Sigma$ from $\Sigma$ in terms of homogeneous 
coordinates $(x_i)_{1\le i\le n}$ associated to the $v_i$ \cite{Cox1992}. First, if $M=N^*$ is the dual lattice and 
$\Pic(\bbP_\Sigma) \cong \bbZ^{k}$ the Picard lattice, the exact sequence
\begin{equation}
 \label{exactsequence}
 M \stackrel{v^*}{\rightarrow}  
 \bbZ^{n} \stackrel{w}{\rightarrow} \Pic(\bbP_\Sigma)
\end{equation}
induces a $\bbZ^k$-grading on the homogeneous coordinate ring through the $(\bbC^\times)^k$-action
$(\bbC^\times)^{k}\ni \lambda:(x_1,\ldots,x_n) \mapsto (\lambda^{w_1} x_1,\ldots,\lambda^{w_n}x_n)$.
Second, if by a (Batyrev) primitive collection \cite{BC1994} we mean a collection of generators $v_i$ that generates 
none of the cones in $\Sigma$, whereas any proper subset of it does, we define the exceptional set as the
union
$$
Z_\Sigma = \bigcup_{p} Z_{\cali_p}, 
$$
where $p$ indexes all primitive collections, and $Z_{\cali_p} = \cap_{i|v_i\in\cali_p}\{x_i = 0\}$.
Then, the toric variety is the quotient
$$
\bbP_\Sigma = \frac{\bbC^n - Z_\Sigma}{(\bbC^\times)^{k}}.
$$

Given $\bbP_\Sigma$, we consider complete intersections $X \subset \bbP_\Sigma$ defined by transversal
polynomials $G_a$ of $\bbZ^k$-degree $d_a$ for $a=1,\ldots,\ell$. We require $X$ to be a Calabi--Yau 
variety, i.e. $\sum_{i=1}^n w_i = \sum_{a=1}^\ell d_a$.

Switching to the algebraic description, let $R=\bbC[x_1,\ldots,x_n]$ be the $\bbZ^k$-graded
coordinate ring associated with the toric data, and 
$J_\Sigma= \langle\prod_{i|v_i\not\in \sigma} x_{i} | \sigma \in \Sigma\rangle$ 
the Cox ideal \cite{Cox1992}, whose vanishing locus is $Z_\Sigma = {\rm V}(J_\Sigma)$ \cite{BC1994}. The complete intersection ring is
denoted by $S=R/(G_1,\ldots,G_\ell)R$, and the image of the Cox ideal along $R\rightarrow S$ by
the same symbol, $J_\Sigma$. The following definitions involving $J_\Sigma$ can be made over $R$ or
over $S$.

Let $\grS$ be the abelian category of graded $S$-modules. The morphisms of $\grS$ are the module
homomorphisms of degree $0$. For a given $S$-module $A$, and $q\in \bbZ^k$, a shift in degree is denoted
by $A(q)$. We call an $S$-module $J_\Sigma$-torsion if it is annihilated by  $(J_\Sigma)^m$ for some 
positive integer $m$.

\begin{definition}
  Let $D^b(\grS)$ be the bounded derived category of graded $S$-modules and $D^b(\torS)$ the full 
triangulated subcategory of graded $J_\Sigma$-torsion $S$-modules, i.e. any object is isomorphic to a complex of $J_\Sigma$-torsion modules. 
The quotient category,
\begin{equation}
\label{quotient}
    D^b(\qgrS)=\frac{D^b(\grS)}{D^b(\torS)},
\end{equation}
is defined by localization along the multiplicative system of morphisms $s$ 
that fit into distinguished triangles
  $$
    A \stackrel{s}{\longrightarrow} B \longrightarrow C \longrightarrow A[1],
  $$
where $A,B \in D^b(\grS)$ and $C \in D^b(\torS)$.
\end{definition}

A generalization of Serre's correspondence identifies the abelian category of coherent sheaves on $X \subset \bbP_\Sigma$ with the abelian category of graded $S$-modules modulo $J_\Sigma$-torsion modules. This correspondence
extends to the derived categories, 
$$
D^b(X) \cong D^b(\qgrS).
$$ 
When we subsequently consider the derived category of coherent sheaves on $X$, we will stick to its 
algebraic description in terms of modules over rings, that is, to $D^b(\qgrS)$.

Since any graded $S$-module has a projective resolution 
(though possibly unbounded from below), 
the derived category $D^b(\grS)$ can be defined to be the homotopy category of projective complexes 
$K^-(\grS)$ with quasi-isomorphisms inverted. Here, the minus index indicates that the complexes may 
be unbounded to the left, but with bounded cohomology. Namely, in the homotopy category of graded $S$-modules all quasi-isomorphisms 
are homotopy equivalences and therefore $D^b(\grS)$ is in fact $K^-(\grS)$. 

By taking the quotient (\ref{quotient}), all objects of the full subcategory $D^b(\torS)$ become zero 
objects in $D^b(\qgrS)$. Over $R$, an important class of zero objects is given by the Koszul complexes 
associated with the primitive collections $\cali_p$,
$$
\calk_p(R):=\calk(\{x_i\}_{i|v_i\in\cali_p};R) = \bigotimes_{i|v_i\in\cali_p} \cone(x_i:R \rightarrow R(w_i))\ .
$$
Indeed, the Koszul complex $\calk_p(R)$ is nothing but the resolution of a $J_\Sigma$-torsion $R$-module,
hence $\calk_p(R)\cong 0$ in $D^b(\qgrR)$.
Over the complete intersection ring $S$ an even shorter regular sequence 
$\cals_{\cali_p} \subset \{x_i\}_{i|v_i\in\cali_p}$ may be the resolution of a 
$J_\Sigma$-torsion $S$-module. For $G_1,\ldots, G_\ell$ generic, $\cals_{\cali_p}$ is just a subset of
$\{x_i\}_{i|v_i\in\cali_p}$, with corresponding set of divisor clases $\{D_j\}_j$. In general,
we can choose representatives $s_j$ of that shorter list of divisor classes, such
that $\bigl(\cap_j\{s_j=0\}\bigr)\cap \bigl(\cap_a\{G_a=0\}\bigr)$ has support on $Z_{\cali_p}$.
Thus, we have
\begin{lemma}
\label{CICYKoszul}
Let $S$ be the complete intersection ring associated with the regular sequence $(G_1,\ldots,G_\ell)$. For 
a primitive collection $\cali_p$ of the smooth toric fan $\Sigma$, we denote by $\cals_{\cali_p}$ a regular
sequence as described above. We let $w_j$ be the degree of $s_j$. Then 
the Koszul complex
$$
  \calk(\cals_{\cali_p};S) =
  \bigotimes_{j} 
  \cone(s_j:S \rightarrow S(w_j))\ ,
$$
is a zero object in $D^b(\qgrS)$. \qed
\end{lemma}

The goal of this work is to understand autoequivalences of $D^b(\qgrS)\cong 
D^b(X)$ and relations among them that can be directly traced back to the toric data, and more specifically to the secondary fan
of $X\subset \bbP_\Sigma$. To this end, we need a realization of $D^b(\qgrS)$ 
that has a natural generalization to all other maximal cones in the secondary fan. This 
alternative construction begins with realizing $X$ as the critical locus $\{{\rm d}W=0\}$ of a  
holomorphic function $W$ on the total space of a certain holomorphic vector bundle over $\bbP_\Sigma$. This then 
leads us to consider the singularity category associated with $W$ \cite{Orlov2004,Orlov2005}.

\section{The singularity category for complete intersections}
\label{CICY}

We append the data of the $\ell$ polynomials $G_a$ to the toric data of $\bbP_\Sigma$ to obtain an enhanced fan $\hat\Sigma$ in $\hat N_\bbR$, where $\hat N = N\oplus \bbZ^{\ell}$, cf. \cite[Chapter 5]{CK2000}. Explicitly, if $v_i$ has coordinates $(v_i^1,\ldots,v_{i}^{n-k})$ with respect to some basis of $N$, the generators $\hat v_i$ of the one-dimensional cones of $\hat\Sigma$ are as follows. For $i=1,\ldots,n$ the coordinates of $\hat v_i$ are $(v_i^1,\ldots,v_i^{n-k},u_{i}^1,\ldots,u_{i}^\ell)$, where the $u_i^a$'s are chosen to satisfy $\sum_{i=1}^n w_i u_{i}^{a} = d_a$. For $a=1,\ldots,\ell$, $\hat v_{n+a}$ is given by a vector with $1$ at the $(n+a)$-th position and $0$ else.
Let $\hat w$ subsume the vectors $(w_1,\ldots,w_n,-d_1,\ldots,-d_{\ell})$. We will sometimes use $w_{n+a} := -d_a$.
Then the exact sequence (\ref{exactsequence}) is extended to
$$
  \hat M \stackrel{\hat v^*}{\rightarrow}  
  \bbZ^{n+\ell} \stackrel{\hat w}{\rightarrow} 
  \Pic(\bbP_{\hat \Sigma}) = 
  \Pic(\bbP_\Sigma) \ .
$$
The (Calabi--Yau) toric variety of the enhanced fan $\hat\Sigma$ is given by
$$
  \bbP_{\hat\Sigma} = \frac{\bbC^{n+\ell} - Z_{\hat\Sigma}}{(\bbC^\times)^{k}} \cong
  \mathrm{Tot}\left(\oplus_{a=1}^\ell \calo(-d_a) \rightarrow \bbP_\Sigma \right)\ .
$$
Note that the exceptional set $Z_{\hat\Sigma}$ is just the pull-back of $Z_\Sigma$ along $\bbC^{n+\ell} 
\rightarrow \bbC^n$.

We denote by $p_a$ the homogeneous coordinate associated with $\hat v_{n+a}$ in $\hat \Sigma$. Let 
$\hat R = \bbC[x_1,\ldots,x_n,p_1,\ldots,p_\ell]$ be the $\bbZ^k$-graded homogeneous coordinate ring 
associated with the fan $\hat\Sigma$. Similarly, $J_{\hat\Sigma}$ be the Cox ideal for $\hat\Sigma$. 
For what follows it is necessary to introduce an additional $2\bbZ$-grading, called R-grading, on 
$\hat R$. The R-grading of the $x_i$ is $0$, that of the $p_a$ is $2$. For a given $\hat R$-module $A$, 
a shift in R-grading by $2r$ is denoted by $A\{2r\}$.

The holomorphic function $W$ on $\bbP_{\hat \Sigma}$ is built from the regular sequence 
$(G_1,\ldots,G_\ell)$ and the auxiliary coordinates $p_a$ as
$$
  W = \sum_{a=1}^\ell p_a G_a\ .
$$
$W$ is a polynomial of degree $0$ and R-grading $2$.
Letting $\hat S=\hat R/(W){\hat R}$, we have the isomorphism of graded rings,
$S\cong \hat S/(G_1,p_1,\ldots,G_\ell,p_\ell)\hat S$, which geometrically corresponds to the 
embedding of the complete intersection $X$ in the toric Calabi--Yau variety $\bbP_{\hat\Sigma}$ 
as the critical locus of $W$,
$$
  X 
    = \{{\rm d}W=0\} \subset \bbP_{\hat\Sigma}.
$$
This follows on account of the transversality of the polynomials $G_1,\ldots,G_\ell$, because 
${\rm d}W$ vanishes iff $p_a=0$ 
and $G_a=0$ for $a=1,\ldots,\ell$.

Following \cite{Orlov2004,Orlov2005}, we now introduce the singularity category
associated with the polynomial $W$. 
\begin{definition}
  \label{defSing}
Let $D^b(\grSh)$ be the derived category of graded $\hat S$-modules, and $\Perf(\grSh)$ the full 
triangulated subcategory of perfect complexes, i.e. bounded complexes of free modules. The singularity 
category is the quotient
  $$
    D_{\sing}(\grSh) := \frac{D^b(\grSh)}{\Perf(\grSh)}.
  $$
\end{definition}

By a result due to Eisenbud \cite{eisenbud}, a (minimal) free resolution of any $\hat S$-module becomes 
two-periodic (up to a shift of R-grading) after a finite number of steps. (The same is true for complexes 
of $\hat S$-modules by induction on the length of the complex.) Modding out by perfect complexes means that 
we can cut off and add finite pieces from the infinite resolution. Therefore, any non-zero object in 
$D_\sing(\grSh)$ can be represented by a half-infinite free complex,
\begin{equation}
  \label{Mcomplex}
  A^\cdot = \ldots 
  \stackrel{f_A}{\longrightarrow} A_0\{-2\}
  \stackrel{g_A}{\longrightarrow} A_1\{-2\}
  \stackrel{f_A}{\longrightarrow} A_0
  \stackrel{g_A}{\longrightarrow} \underline{A_1},
\end{equation}
with ${\rm rk} A_0 = {\rm rk} A_1$. The underlined module is in homological degree $0$. 

Notice that 
\begin{equation}
  \label{Rshift}
  A^\cdot\{2\} = A^\cdot[2],
\end{equation}  
for any $A^\cdot \in D_\sing(\grSh)$. The shift in R-grading 
by $-2$ in (\ref{Mcomplex}) is due to the R-grading of $W$. In fact, the 
complex of free $\hat S$-modules (\ref{Mcomplex}) can be lifted to a sequence of free $\hat R$-modules so 
that the homomorphisms compose as
$f_{\tilde A} \, g_{\tilde A} = W \cdot \id_{\tilde A_0}$ and 
$g_{\tilde A} \, f_{\tilde A} = W \cdot \id_{\tilde A_1}$, where 
$A_j = \tilde A_j \otimes_{\hat R} \hat S$ for $j=0,1$.
This also leads to Orlovs result in \cite{Orlov2005} that $D_\sing(\grSh)$ can equivalently be described by the 
triangulated category of matrix factorizations of $W$ over the ring $\hat R$. 

The morphisms in the  singularity category are given by chain maps modulo homotopy and modulo chain maps that 
factor through perfect complexes. On the two-periodic part 
the 
morphisms are also two-periodic, and since we mod out by chain maps that factor through perfect complexes, 
the two-periodic part determines the morphisms uniquely up to 
homotopy. On the 
representatives (\ref{Mcomplex}) a chain map is a commutative diagram
\begin{equation}
  \label{Mmorphism}
  \begin{array}{cccccccc}
  \ldots &A_0\{-2\}&
  \stackrel{g_A}{\longrightarrow} &A_1\{-2\}&
  \stackrel{f_A}{\longrightarrow} &A_0&
  \stackrel{g_A}{\longrightarrow} &\underline{A_1} \\[5pt]
  & ~~\downarrow \psi_0 && ~~\downarrow \psi_1 &
  & ~~\downarrow \psi_0 && ~~\downarrow \psi_1\\
  \ldots &B_0\{-2\}&
  \stackrel{g_B}{\longrightarrow} &B_1\{-2\}&
  \stackrel{f_B}{\longrightarrow} &B_0&
  \stackrel{g_B}{\longrightarrow} &\underline{B_1}
  \end{array}
\end{equation}

For the following purposes, a useful representative of the object $A^\cdot \in D_\sing(\grSh)$ in $D^b(\grSh)$ 
is constructed by continuing the 
complex (\ref{Mcomplex}) periodically by $2r$ steps to the right, so that the R-degrees of $A_0\{2r\}$ 
and $A_1\{2r\}$ are all positive, and then cutting off the finite piece with positive R-degrees, in any 
homological degree. We denote the corresponding functor by
$$
  \sigma_{\leq 0} : D_\sing(\grSh)\longrightarrow D^b(\grSh).
$$
Let $\Psing:D^b(\grSh) \longrightarrow D_\sing(\grSh)$, then clearly we have $A^\cdot \cong \Psing 
\sigma_{\leq 0}A^\cdot$ for any object $A^\cdot \in D_\sing(\grSh)$. It follows from two-periodicity 
that a free resolution in the image of $\sigma_{\leq 0}$, say $B^\cdot$, satisfies the relation 
\begin{equation}
  \label{standardform}
  \sigma_{\leq 0} \Psing(B^\cdot\{2r\})= B^\cdot[2r], \qquad
  \mathrm{for}\quad r=0,1,2,\ldots\ .
\end{equation}
In fact, the subcategory of such objects is a full triangulated subcategory of $D^b(\grSh)$.

\begin{lemma}
  \label{equiv standard form}
  Let $D_{\leq 0}$ be the full triangulated subcategory of objects satisfying
  (\ref{standardform}). Then the adjoint pair of functors,
  $$
    D_\sing(\grSh) \quad\mapup{\sigma_{\leq 0}}\mapback{\Psing}\quad 
    D_{\leq 0} ~~\subset~~ D^b(\grSh),
  $$
  is an equivalence of triangulated categories.
\end{lemma}
\proof We know already that $A^\cdot \cong \Psing \sigma_{\leq 0}A^\cdot$ for any object $A^\cdot \in 
D_\sing(\grSh)$, and by definiton $B^\cdot \cong \sigma_{\leq 0}\Psing B^\cdot$ for any $B^\cdot \in 
D_{\leq 0}$. It remains to check isomorphism of morphisms: Because of two-periodicity 
a chain map between objects in $D_{\leq 0}$ cannot factor through a perfect complex, and hence
$$
  \Hom{}_{\sing}(\Psing A^\cdot,\Psing B^\cdot) ~~\cong  ~~
  \Hom{}_{D_{\leq 0}}(A^\cdot,B^\cdot),
$$
for $A^\cdot,B^\cdot \in D_{\leq 0}$. \qed

\begin{proposition}
\label{DtoMF}
The triangulated categories $D_\sing(\grSh)$ and $D^b(\grS)$ are equivalent.
\end{proposition}

\proof
Using Lemma~\ref{equiv standard form} it remains to construct an adjoint pair of functors,
$$
\begin{array}{rcccl}
  D^b(\grR) ~~\supset\quad D^b(\grS) \quad
   \mapup{E}\mapback{\omega}\quad  D_{\leq 0} \quad\subset~~ D^b(\grSh)
  \end{array}
$$\\
that implements the equivalence. 

The functor $E$ is provided by Eisenbud's work \cite{eisenbud} and is constructed in two steps. 
For the first step we include $D^b(\grS)$ as full triangulated subcategory in $D^b(\grR)$ by considering 
$S$-modules as $R$-modules which are annihilated by the regular sequence $(G_1,\ldots,G_\ell)$. For an 
object $A\in D^b(\grS)$ we take the (minimal) $R$-free resolution $P^\cdot(A)$. Let $P(A)$ be the free 
$R$-module in the complex $P^\cdot(A)$ and $d_0$ the differential.

The second step constructs from $P^\cdot(A)$ a half-infinite complex in $D_{\leq 0}$. Theorem 7.1 of 
\cite{eisenbud} allows us to introduce auxiliary endomorphisms $d_{\bfn}:P(A)[2|\bfn|] \rightarrow P(A)$ 
for $\bfn \in (\bbZ_{\geq 0})^\ell$ and $|\bfn|= \sum_{a=1}^\ell n_a$, which have homological degree $1$ 
as well as degree $\sum_a n_a d_a$ and satisfy
\begin{equation}
    \label{recursiveMF}
    \begin{array}{rccl}
    d_{\bf 0}d_{\bf e_a}+d_{\bf e_a}d_{\bf 0} &=& G_a \cdot \id_{P(A)}[-2], \\
    \sum_{{\bf m},|{\bf m}|\leq|\bfn|} d_{{\bf m}}d_{\bfn-{\bf m}}&=& 0,&
    \qquad\mathrm{for}\quad |\bfn|>1.
    \end{array}
\end{equation}
Here $e_a = (0,\ldots,0,1,0,\ldots,0)$ with the $1$ at the a${}^\mathrm{th}$
position, and we set $d_{\bf 0} = d_0$. 

Define the free $\hat S$-module 
$\hat P(A) = (\bbC[p_1,\ldots,p_\ell] \otimes_\bbC P(A)) / (W)$. Then 
$E(A) \in D_{\leq 0}$ is the $\hat S$-module
\begin{equation}
  \nonumber
  \oplus_{r=0}^\infty \hat P(A)\{-2r\}[2r]
\end{equation}
together with the differential $d = \sum_\bfn p^\bfn \otimes d_\bfn$. Notice that the latter has degree $0$ and 
preserves the R-grading. In view of (\ref{recursiveMF}), $d^2=W\cdot \id [-2]$, which is zero on $\hat 
S$-modules.
$E(A)$ satisfies condition (\ref{standardform}) and is therefore an object in 
$D_{\leq 0}$. In fact, the object is determined by $d_0$ uniquely up to isomorphism, since finding a 
solution to the recursive relations (\ref{recursiveMF}) admits the freedom, $d \rightarrow U d U^{-1}$ for 
$U = \id + \sum_{\bfn > 0} p^\bfn u_\bfn$.

The functor $\omega$ is defined by
\begin{eqnarray}
  \nonumber
  \omega : D_{\leq 0} &\rightarrow & D^b(\grS) \subset D^b(\grR), \\
  \nonumber
  A^\cdot\quad &\mapsto& (\sigma_{\geq 0} A^\cdot) \otimes 
  \hat S/(p_1,\ldots,p_\ell)\hat S.
\end{eqnarray}
$\sigma_{\geq 0}$ cuts off negative R-gradings of $A^\cdot$ and therefore picks the R-grading $0$ component, 
whose differential does not contain the auxiliary coordinates $p_a$. Tensoring with $\hat S/(p_a)\hat S 
\cong R$ then removes the auxiliary coordinates from the module as well, i.e. the image of the functor 
$\omega$ is in $D^b(\grR)$. Indeed, the image of $\omega$ is $D^b(\grS)$: Given an object $A^\cdot \in 
D_{\leq 0}$ we write its differential as $d = \sum_\bfn p^\bfn d_\bfn$. Then Theorem 7.2 of \cite{eisenbud} 
constructs from the endomorphisms $d_\bfn$ an infinite $S$-free resolution, that is an object in $D^b(\grS)$, 
which is isomorphic to $\omega(A^\cdot)$ in $D^b(\grR)$. 

From the definitions of $E$ and $\omega$ we have 
$$
  \omega E (A) = P^\cdot(A) \cong A,
$$
for any $A \in D^b(\grS) \subset D^b(\grR)$.
Conversely, for $B^\cdot\in D_{\le 0}$, since $\omega(B^\cdot)$ is the R-grading $0$ component of 
$B^\cdot$ and $E$ reconstructs the auxiliary endomorphisms in $B^\cdot$, we also have 
$$
  E\omega (B^\cdot) \cong B^\cdot.
$$ 

It remains to show that $\omega$ and $E$ are adjoint functors, that is
\begin{equation*}
  \Hom{}_{D_{\leq 0}}(E(A),B^\cdot) =
  \Hom{}_{\grS}(A,\omega(B^\cdot)),
\end{equation*}
for $A \in D^b(\grS)$ and $B^\cdot \in D_{\le 0}$. We use that $D^b(\grS)$ is full in $D^b(\grR)$ 
and write the right-hand side as $\Hom{}_{\grR}(P^\cdot(A),\omega(B^\cdot))$. For the left-hand side 
write a morphism in $\Hom{}_{D_{\leq 0}}(E(A),B^\cdot)$ as $\psi = \sum_{\bfn\in\bbZ_{\geq 0}} p^\bfn 
\psi_\bfn$. Just as the differential $d$, the morphism $\psi$ is determined (up to isomorphisms) 
by its R-grading $0$ component $\psi_0$, hence the left-hand side is also isomorphic to 
$\Hom{}_{\grR}(P^\cdot(A),\omega(B^\cdot))$. \qed

\begin{definition}
Let $D^b(\torSh)$ be the full triangulated subcategory of $D^b(\grSh)$ consisting of graded 
$J_{\hat \Sigma}$-torsion modules, and consider its image  in $D_\sing(\grSh)$, that is 
$D_\sing(\torSh):=\Psing D^b(\torSh)$. This is a full triangulated subcategory of $D_\sing(\grSh)$. We 
define the quotient category
$$
D_\sing(\qgrSh) := \frac{D_\sing(\grSh)}{D_\sing(\torSh)}.
$$
\end{definition}

\begin{theorem}
\label{DtoMFquotient}
The image of the torsion subcategory $D_\sing(\torSh)$ under the equivalence 
$\omega_{\le 0}=\omega\sigma_{\le 0}: D_\sing(\grSh) \rightarrow D^b(\grS)$ is $D^b(\torS)$, 
and, consequently,
\begin{equation}\nonumber
 D_\sing(\qgrSh) ~\cong~ D^b(\qgrS).
\end{equation}
\end{theorem}

\proof Recalling that $Z_{\hat\Sigma}$ is the pull-back of $Z_\Sigma$ we find from the definition of $\omega$ that it maps $J_{\hat \Sigma}$-torsion complexes to 
$J_{\Sigma}$-torsion complexes, hence $\omega_{\le 0}(D_\sing(\torSh)) \subseteq D^b(\torS)$. Also, 
from the definition of $D_\sing(\torSh)$, we know that $\Psing (D^b(\torS))\subseteq D_\sing(\torSh)$. 
Applying $\omega_{\le 0}$, we obtain $\omega_{\le 0}\circ \Psing (D^b(\torS)) = D^b(\torS) \subseteq 
\omega_{\le 0} (D_\sing(\torSh))$, and therefore $D^b(\torS) \cong \omega_{\le 0} (D_\sing(\torSh))$. 
\qed

\section{Phases of triangulated categories}
\label{phases}

The realization of the derived category $D^b(X)\cong D^b(\qgrS)$ 
through the singularity category $D_\sing(\qgrSh)$ in Theorem~\ref{DtoMFquotient} 
motivates us to define singularity categories for each maximal cone of the secondary 
fan of $\bbP_{\hat\Sigma}$. We shall label the latter by the toric fan $\hat\Phi$.

The corresponding toric variety is 
\begin{equation}
  \label{toricVar}
  \bbP_{\hat\Phi} = \frac{\bbC^{n+\ell}-Z_{\hat\Phi}}{(\bbC^\times)^k}.
\end{equation}
Note that in general $\hat\Phi$ contains fewer one-dimensional cones than $\hat\Sigma$, i.e. fewer than 
$n+\ell$. This is taken into account in (\ref{toricVar}) by the exceptional set being the union, 
$$
  Z_{\hat\Phi} = Z_{\hat\Phi}^{\it prim} \cup Z_{\hat\Sigma(1)\setminus\hat\Phi(1)}.
$$
Here $Z_{\hat\Phi}^{\it prim}$ is given in terms of Batyrev's primitive collections for the fan $\hat\Phi$, 
and\footnote{Note that the union may include some of the auxiliary coordinates $x_{n+a}=p_a$.}
$$
Z_{\hat\Sigma(1)\backslash\hat\Phi(1)} = 
\bigcup_{\hat v_i \in \hat\Sigma(1)\backslash\hat\Phi(1)}\{x_i = 0\}.
$$
Note that
$$
  \bbP_{\hat\Phi} = \frac{(\bbC^{n+\ell-r}-Z_{\hat\Phi}^{prim})\times 
  (\bbC^\times)^r}{(\bbC^\times)^k}
  \cong \frac{(\bbC^{n+\ell-r}-Z_{\hat\Phi}^{prim})}{(\bbC^\times)^{k-r}\times
  G_{\hat\Phi}}
$$
where $r$ is the number of one-dimensional cones in $\hat\Sigma(1)
\backslash\hat\Phi(1)$ and $G_{\hat\Phi}$ is a finite group.

To construct quotient singularity categories in every maximal cone $\hat\Phi$, it is
convenient to start with the graded coordinate rings $\hat R=\bbC[x_1,\ldots,x_n,p_1,\ldots,p_\ell]$ 
and $\hat S = \hat R/(W)\hat R$, with the degree and R-grading as before. In a maximal cone $\hat\Phi$ 
of the secondary fan consider the ideal $J_{\hat\Phi} =
\langle \prod_{i|v_i \in \hat\Sigma(1),v_i\not\in \sigma} x_i|\sigma \in \hat\Phi\rangle$. 
Its vanishing locus is the exceptional set $Z_{\hat\Phi}$. We say that a graded $\hat S$-module is $J_{\hat\Phi}$-torsion if it is annihilated 
by $J_{\hat\Phi}^{\otimes m}$ for some positive integer $m$.

\begin{definition}
  \label{singcat}
  Let $\hat S = \hat R/(W)\hat R$, and $D^b(\grSh)$ as well as 
  $D_\sing(\grSh)$ be as in Definition~\ref{defSing}. Take any 
  maximal cone $\hat\Phi$ in the secondary fan of
  $\bbP_{\hat \Sigma}$. Let $D^b(\torPhi)$ be the full 
  triangulated subcategory of $J_{\hat\Phi}$-torsion modules 
  in $D^b(\grSh)$, and $D_\sing(\torPhi)$ its image in 
  $D_\sing(\grSh)$. Then we define the quotient
  $$
   \calc_{\hat\Phi}= D_\sing(\qgrPhi) := \frac{D_\sing(\grSh)}{D_\sing(\torPhi)}.
  $$
\end{definition}

\begin{remark}
  \label{hatSKoszul}
  For every irreducible component $Z_{\cali_p}$ of the exceptional 
  set $Z_{\hat\Phi}$, the associated Koszul complex 
  $$
    \calk_p(\hat S):=\calk(\{x_i\}_{v_i\in\cali_p};\hat S) = 
    \bigotimes_{i\in\cali_p} \cone(x_i:\hat S \rightarrow \hat S(w_i)\{r_i\})\ .
  $$ 
  is an object in $D^b(\torPhi)$. Here, $r_i$ is the R-degree of $x_i$,
  i.e. $r_i= 0$ for $i=1,\ldots,n$ and $r_{n+a}=2$ for $a=1,\ldots,\ell$.
  Furthermore, for any object $A^\cdot$ 
  of $D^b(\grSh)$, the tensor product $A^\cdot \otimes_{\hat S} \calk_p$
  is in $D^b(\torPhi)$ and via $\Psing$ in $D_\sing(\torPhi)$.
\end{remark}

We have introduced the quotient singularity categories separately for each maximal cone. The following 
result relates them.

\begin{theorem} 
 \label{allphases}
For any pair of neighboring maximal cones, $\hat\Phi_1$ and $\hat\Phi_2$, in the secondary fan of $\bbP_{\hat\Sigma}$, with Calabi--Yau complete intersecion $X\subset \bbP_{\hat\Sigma}$, there is a family $\{F_m^{\hat\Phi_2\hat\Phi_1}\}_{m\in\bbZ}$ of 
equivalences of the corresponding quotient singularity categories,
  \begin{equation}
    \label{allphasesequivalence}
    F_m^{\hat\Phi_2\hat\Phi_1}:
   D_\sing({\rm qgr}_{\hat\Phi_1}\hat S)  
   \stackrel{\cong}{\longrightarrow}
   D_\sing({\rm qgr}_{\hat\Phi_2}\hat S).
  \end{equation}
\end{theorem}

\proof  Let $\bbP_{\hat\Phi_1}$ and $\bbP_{\hat\Phi_2}$ be the (Calabi--Yau) toric 
  varieties in the neighboring maximal cones. We let $T\in\Pic(\bbP_{\hat\Sigma})^*$
be the primitive dual vector characterizing the face between the two adjacent maximal cones by the condition $T(v)=0$ for all $v$ in the face.
The exceptional sets $Z_1$ and $Z_2$ in the geometric construction of 
$\bbP_{\hat\Phi_1}$ and $\bbP_{\hat\Phi_2}$ are related by \cite[Chapter 4.5]{HHP2008}
$$
\begin{array}{rcl}
Z_1 &=& Z_+ \cup (Z_1\cap Z_2), \\
Z_2 &=& Z_- \cup (Z_1\cap Z_2),
\end{array}
$$
where $Z_{\pm} = \cap_{i|T(w_i)\gtrless 0 } \{ x_i = 0 \}$.
We let ${\rm tor}_{\hat\Phi_+}\hat S$, ${\rm tor}_{\hat\Phi_-}\hat S$, and 
${\rm tor}_{\hat\Phi_{12}}\hat S$ denote the categories of torsion modules associated 
with $Z_+$, $Z_-$, and $Z_1\cap Z_2$, respectively. 

By the Calabi--Yau condition, $\sum_{i=1}^{n+\ell} w_i =0$, we may define
$$
  \sigma ~:= \sum_{i|T(w_i)>0} \!\!T(w_i) ~~=~~ - \sum_{i|T(w_i)<0} \!\!T(w_i).
$$
For $m \in \bbZ$ let $K^m({\rm gr}\hat S) \subset K^-({\rm gr}\hat S)$ be the homotopy category generated from invertible modules $\hat S(q)$ satisfying (cf. \cite{HHP2008,vdB2002,kawamata})
\begin{equation}
  \label{GRR}
  m \leq T(q) < m+\sigma.
\end{equation}
Denote the associated singularity category by $D^m_\sing(\grSh)$ and consider
\begin{equation}
  \label{phaseequiv}
  D_\sing({\rm qgr}_{\hat\Phi_+}\hat S) 
  \mapup{\omega_+}\mapback{\pi_+}
  D^m_\sing(\grSh)
  \mapup{\pi_-}\mapback{{\omega_-}}
  D_\sing({\rm qgr}_{\hat\Phi_-}\hat S),
\end{equation}
where $D_\sing({\rm qgr}_{\hat\Phi_\pm}\hat S)$ are the quotients of 
$D_\sing(\grSh)$ by $D_\sing({\rm tor}_{\hat \Phi_\pm}\hat S)$. The functors $\pi_\pm$ are the 
projections of $D_\sing(\grSh)$ to the quotient categories applied to the subcategory $D^m_\sing(\grSh)$. 
Comparing (\ref{GRR}) with $Z_\pm$, we find that the objects in $D^m_\sing(\grSh)$ can not be torsion. 
Also, the objects of $D^m_\sing(\grSh)$ generate $D_\sing({\rm qgr}_{\hat\Phi_\pm}\hat S)$, so that 
$\pi_\pm$ are bijective on the set of objects. $D^m_\sing(\grSh)$ containing no torsion objects also 
means that there are no non-trivial extensions, and the functors $\pi_\pm$ are fully faithful, hence 
equivalences. The inverses $\omega_\pm$ take the unique (up to homotopy) representative of an 
isomorphism class in $D_\sing({\rm qgr}_{\hat\Phi_\pm}\hat S)$ that satisfies the restriction 
(\ref{GRR}).

Finally, taking at each step in (\ref{phaseequiv}) the quotient by 
$D_\sing({\rm tor}_{\hat\Phi_{12}}\hat S)$ we find that the functor in the 
theorem is given by the compositions,
$$
  D_\sing({\rm qgr}_{\hat\Phi_1}\hat S) 
  \stackrel{\omega_+}{\longrightarrow}
  D^m_\sing({\rm qgr}_{\hat\Phi_{12}}\hat S)
  \stackrel{\pi_-}{\longrightarrow}
  D_\sing({\rm qgr}_{\hat\Phi_2}\hat S).
$$
\qed

  A combination of Theorems~\ref{DtoMFquotient} and \ref{allphases} 
  was proved for the hypersurface case by Orlov 
  in \cite{Orlov2005}. Van den Bergh \cite{vdB2002} and Kawamata 
  \cite{kawamata} stated analogous results for the derived categories
  $D^b(\bbP_{\hat\Phi})$ of 
  Calabi--Yau toric varieties using \emph{noncommutative} crepant resolutions. 
  They build a tilting module, say $P$, out of a generating set of invertible 
  modules $\hat S(q)$ satisfying (\ref{GRR}) and use the derived category of 
  the noncommutative endomorphism algebra, $D^b(\End(P))$, to show the 
  equivalence of $D^b(\bbP_{\hat\Phi_1})$ and $D^b(\bbP_{\hat\Phi_2})$.

\begin{remark}
\label{stringy}
  In string theory the objects of the quotient singularity categories correspond 
  to boundary condition of certain two-dimensional supersymmetric field theories. 
  A physics derivation of the functor (\ref{allphasesequivalence})
  was given in \cite{HHP2008}. Therein it was found that the choice of the integer $m$ 
  corresponds to a choice of a homotopy class of paths connecting limit points 
  corresponding to $\hat\Phi_1$ and $\hat\Phi_2$ in the moduli space $B$ of the 
  mirror of $X$.
\end{remark}

\section{Autoequivalences: The main formula}
\label{main}

\def\cala{{\mathcal A}}

In this section we study elements of the group $\aut_{\hat\Phi}$ of autoequivalences of 
$D_\sing(\qgrPhi)$. We begin in the phase $\hat\Phi=\hat\Sigma$.

Immediate elements of $\aut_{\hat\Sigma}$, the group of automorphisms of $D^b(X)
\cong D^b(\qgrS)$, are given by the shift functors $[n]$, for $n \in \bbZ$, and the twists $\M^q$ by 
the modules $S(q)$ for any $q \in \bbZ^k \cong \Pic(\bbP_{\hat\Sigma})$. For $A \in D^b(\qgrS)$, 
$$
  \M^q : A \mapsto A \otimes_S S(q).
$$
Together with the automorphisms $\aut(X)$ of $X$, that is, the graded ring automorphisms of $S$ 
modulo $(\bbC^\times)^k$, these generate a subgroup 
$\cala_0(X) \cong \aut(X) \ltimes \Pic(\bbP_{\hat\Sigma}) \times \bbZ$ of $\aut(D^b(X))$.

Because of Theorem~\ref{DtoMFquotient} the group of automorphisms $\aut(D^b(X))$, and in particular 
its subgroup $\cala_0(X)$, also acts on the quotient singularity category $D_\sing(\qgrSh)$. The latter category appears to have an 
additional functor of twisting by $\hat S\{2r\}$. Recall however from (\ref{Rshift}) that a shift 
in R-grading by $\{2r\}$ is equal to the shift functor $[2r]$, hence does not introduce a new 
autoequivalence.

Let $s$ be an element of $S$ with degree $w$. Introduce on $D^b(\qgrS)$ the endofunctor
$$
  \N(s): A \mapsto \cone\left(s:A \rightarrow \M^w (A)\right).
$$

\begin{proposition} 
  \label{twistmonoRelation}
   For every primitive collection $\cali_p$ of $\Sigma$ choose a 
   regular sequence $\cals_{\cali_p}$ of elements $s_j \in S$, 
   as in Lemma~\ref{CICYKoszul}. Then,
	 \begin{equation}
	   \label{0relation}
	   \bigcirc_{s_j \in \cals_{\cali_p}} \N(s_j)(A) \cong 0 , 
   \end{equation}
   for any $A\in D^b(\qgrS)$.
\end{proposition}

\proof By Lemma~\ref{CICYKoszul} the Koszul complex $\calk(\cals_{\cali_p};S)$ is isomorphic 
to the zero object in $D^b(\qgrS)$. The same is true for the tensor product $A \otimes_S 
\calk(\cals_{\cali_p};S)$ for any object $A$. This is nothing but the left-hand side of 
(\ref{0relation}) since $\N(s)(A) = A \otimes_S\cone\left(s:S \rightarrow S(w)\right)$. 
\qed

For any other maximal cone $\hat\Phi$ of the secondary fan, there are twist autoequivalences acting on $D_\sing(\qgrPhi)$ as well. We set
$$
  \M_{\hat\Phi}^q: A^\cdot \mapsto A^\cdot \otimes_{\hat S} \hat S(q)
  \qquad \mathrm{for}\quad A^\cdot \in D_\sing(\qgrPhi).
$$
For any element $s \in \hat S$ with degree $w$ and R-grading $2r$, let
$$
  \N_{\hat\Phi}(s): A^\cdot \mapsto
  \cone\left(s:A^\cdot \rightarrow \M_{\hat\Phi}^w (A^\cdot)[2r]\right).
$$
Notice the shift in homological degree and recall that $[2r]\cong\{2r\}$ on $D_\sing(\qgrPhi)$. Then, we have
\begin{proposition}
   \label{twistmonoRelationPhi}
   Let $\hat\Phi$ be the fan associated to a maximal cone in the secondary 
   fan of $\bbP_{\hat\Sigma}$.
   Consider an arbitrary object $A^\cdot$ in the quotient singularity category
   $D_\sing(\qgrPhi)$. 
   Then, for every primitive collection $\cali_p$ of $\hat\Phi$,
   $$
     \bigcirc_{i|\hat v_i \in \cali_p} \N_{\hat\Phi}(x_i)(A^\cdot) \cong 0 , 
   $$
   and for every $\hat v_i \in \hat\Sigma(1)\backslash \hat\Phi(1)$ 
   associated with $x_i \in \hat S$, 
   $$
     \N_{\hat\Phi}(x_i)(A^\cdot) \cong 0,\quad
     \mathrm{or~equivalently,}\quad
     \M_{{\hat\Phi}}^{w_i}(A^\cdot) \cong (A^\cdot)[-2r_i].
   $$
\end{proposition}

\proof Using Remark~\ref{hatSKoszul}, the proof is similar to the proof of Proposition~\ref{twistmonoRelation}.
\qed

\begin{remark}
In the maximal cone $\hat\Sigma$ we may write $\M^q = \omega \sigma_{\leq 0} \M_{\hat\Sigma}^q \Psing E$.
Using this equivalence, the relations of Proposition~\ref{twistmonoRelation} are in general stronger than 
those of Proposition~\ref{twistmonoRelationPhi}. One may sometimes obtain similarly strong relations also 
in the other phases, namely whenever the exceptional set $Z_{\hat\Phi}$ admits the solution of ${\rm d}W=0$ 
to set $p_a=G_a=0$ for (at least) one $a=1,\ldots,\ell$. Then, we may work over the ring 
$\hat S/(p_a,G_a)\hat S$. In fact, Proposition~\ref{twistmonoRelation} applies because over 
$\bbP_{\hat\Sigma}$, ${\rm d}W=0$ admits $p_a=G_a=0$ for \emph{all} $a=1,\ldots,\ell$. 
\end{remark}

Theorem~\ref{allphases} relates the singularity categories 
in neighbouring maximal cones $\hat\Phi_1$ and $\hat\Phi_2$. The 
associated twists are only partly independent autoequivalences. In fact, condition (\ref{GRR}), which 
defines the functor $F_m^{\hat\Phi_2\hat\Phi_1}:\calc_{\hat\Phi_1}\rightarrow \calc_{\hat\Phi_2}$, 
implies that
\begin{equation}
  \label{equalautos}
  \M^q_{\hat\Phi_2} \circ F_m^{\hat\Phi_2\hat\Phi_1}  ~~\cong~~ 
  F_{m+T(q)}^{\hat\Phi_2\hat\Phi_1} \circ \M^q_{\hat\Phi_1} .
\end{equation}
In particular, the twists with $T(q)=0$ commute with the functor $F_m^{\hat\Phi_2\hat\Phi_1}$ (for 
a fixed integer $m$), which says that they correspond to the same autoequivalence,
$$
  \M^q_{\hat\Phi_1} ~~\cong~~ \M^q_{F_m^{\hat\Phi_2\hat\Phi_1}}
  := (F_m^{\hat\Phi_2\hat\Phi_1})^{-1} \M^q_{\hat\Phi_2} 
  \circ F_m^{\hat\Phi_2\hat\Phi_1}, \quad
  \mathrm{for}~~T(q)=0.
$$
On the other hand, if $T(q)$ does not vanish the composition,
\begin{equation}
  \label{STtwist}
  (F_m^{\hat\Phi_2\hat\Phi_1})^{-1}
  F_{m+T(q)}^{\hat\Phi_2\hat\Phi_1} \cong
  \M^q_{F_m^{\hat\Phi_2\hat\Phi_1}} 
  \M^{-q}_{\hat\Phi_1},
\end{equation}
is a non-trivial autoequivalence of $\calc_{\hat\Phi_1}$.

For any $\hat\Phi$ let $F_{\hat\Phi}:\calc_{\hat\Sigma}\rightarrow \calc_{\hat\Phi}$ be a composition 
of functors of Theorem~\ref{allphases}. Then the twists in $\hat\Phi$ act on $\calc_{\hat\Sigma} 
\cong D^b(X)$ via $\M^q_{F_{\hat\Phi}} = F_{\hat\Phi}^{-1} \circ \M_{\hat\Phi}^q \circ F_{\hat\Phi}$.
By combining Proposition \ref{twistmonoRelationPhi} with Theorem \ref{allphases} we obtain
\begin{theorem}
\label{evenbigger}
  Let $X$ be a smooth Calabi--Yau complete intersection in a toric variety
  $\bbP_\Sigma$, and $D_\sing(\qgrSh)$ ($\cong D^b(\qgrS)$) 
  its quotient singularity category. 
  Then for every maximal cone $\hat\Phi$ in the secondary fan of 
  $\bbP_{\hat\Sigma}$ the autoequivalences
  $\M_{F_{\hat\Phi}}^q$ induce an action of the Picard lattice $\Pic(\bbP_{\hat\Sigma})$ 
  on $D_\sing(\qgrSh)$, subject to the following relations. For any object 
$A^\cdot \in D_\sing(\qgrSh)$ and every primitive collection $\cali_p$ of $\hat\Phi$,
   \begin{equation}
    \label{primitiverelation}
     \bigcirc_{i|\hat v_i \in \cali_p} 
     \N_{F_{\hat\Phi}}(F_{\hat\Phi}(x_i))(A^\cdot) \cong 0 , 
   \end{equation}
   and for every $\hat v_i \in \hat\Sigma(1)\backslash \hat\Phi(1)$ 
   and associated element $x_i \in \hat S$, 
   \begin{equation}
     \label{equivalencerelation}
     \N_{F_{\hat\Phi}}(F_{\hat\Phi}(x_i))(A^\cdot) \cong 0,\quad
     \mathrm{or~equivalently,}\quad
     \M_{F_{\hat\Phi}}^{w_i}(A^\cdot) \cong (A^\cdot)[-2r_i].
   \end{equation}
\end{theorem}

\proof Writing for any $A^\cdot \in \calc_{\hat\Sigma}$
\begin{eqnarray}
  \nonumber
  \N_{F_{\hat\Phi}}(F_{\hat\Phi}(s))(A^\cdot) &:=&
  \cone(F_{\hat\Phi}(s):A^\cdot \longrightarrow 
  \M^w_{F_{\hat\Phi}}(A^\cdot)[2r]) = \\
  \nonumber
  &=& F_{\hat\Phi} \cone(s:F_{\hat\Phi}^{-1}(A^\cdot) \longrightarrow
  \M^w_{\hat\Phi} F_{\hat\Phi}^{-1}(A^\cdot)) = \\
  \nonumber 
  &=&
  F_{\hat\Phi} \N_{\hat\Phi}(s) F_{\hat\Phi}^{-1}(A^\cdot),
\end{eqnarray}
the claim follows immediately from Proposition~\ref{twistmonoRelationPhi}.
\qed

\section{Applications}
\label{applications}

\subsection{Monodromy representations --- more about $\aut(D^b(X))$.}\

According to Theorem \ref{evenbigger}, we have an action of the toric Picard lattice of the
complete intersection $X$ on $D^b(X)$ for every phase $\hat\Phi$ in the secondary fan
of $\bbP_{\hat\Sigma}$, generalizing the simple twisting by line bundles for $\hat\Phi=\hat\Sigma$. As is 
evident from the relations (\ref{primitiverelation}) and (\ref{equivalencerelation}), this 
action can be very different for each $\hat\Phi\neq\hat\Sigma$. 

In terms of the monodromy
representation $\pi_1(B)\to\aut(D^b(X))$, mentioned in the introduction,
the action of $\Pic(\bbP_{\hat\Sigma})$ obtained from $\hat\Phi$ is a categorical lift of the monodromies 
around the boundary divisor of $B$ corresponding to the maximal cone $\hat\Phi$ of the secondary 
fan. This point of view fits in nicely with the interpretation of the equivalence between
the various $\calc_{\hat\Phi}$'s as depending on the homotopy class of a path in $B$. According to 
Remark~\ref{stringy} the composition (\ref{STtwist}) for $T(q)=1$ is the monodromy along a path around 
the discriminant locus of $B$, and therefore corresponds to a Seidel--Thomas twist \cite{ST2000}, 
cf. also \cite{Horja2001,AHK2002} and \cite[Chapter 10.5]{HHP2008}.

\subsection{Proof of K-theory formula}\

The projection to K-theory splits the cones in Proposition \ref{twistmonoRelationPhi}, which implies
Corollary \ref{Ktheory}.

When interpreted in terms of monodromy representations (see previous subsection), the relations in K-theory
can also be obtained by studying the analytic continuation of periods of the mirror variety, as in 
\cite{paul}. (See also \cite{AP2009} for a recent study.) Our results provide the precise categorical 
lift of the monodromies and the relations between them.

\subsection{Some examples}\

Consider the list of $14$ Calabi--Yau complete intersections 
$X(d_1\ldots d_\ell)$ in weighted projective spaces $\bbP^{3+\ell}(w_0\ldots w_{3+\ell})$ with rank $1$ 
Picard lattice, that is $\Pic(X) \cong \bbZ$ (cf. \cite{DM2005,CYYE2006}):
\begin{equation}
  \nonumber
  \begin{array}{|rl|rl|rl|}
  \hline
  X(5) \subset& \bbP^4(11111) & X(2,4) \subset& \bbP^5(111111) & X(2,12)\subset& \bbP^5(111146) \\
  X(6) \subset& \bbP^4(11112) & X(3,3) \subset& \bbP^5(111111) & X(4,6) \subset& \bbP^5(111223) \\
  X(8) \subset& \bbP^4(11114) & X(3,4) \subset& \bbP^5(111112) & X(6,6) \subset& \bbP^5(112233) \\
  X(10) \subset& \bbP^4(11125) & X(2,6) \subset& \bbP^5(111113) & X(2,2,3) \subset& \bbP^6(1111111) \\
  &&X(4,4) \subset& \bbP^5(111122) & X(2,2,2,2) \subset& \bbP^7(11111111)
  \\
  \hline
  \end{array}
\end{equation}\\
The polynomial ring associated to the toric data is 
$\hat R = \bbC[x_0,\ldots,x_{3+\ell},p_1,\ldots,p_\ell]$. The secondary fan has two maximal cones, 
say $\hat\Sigma$ and $\hat\Xi$, where the first shall correspond to the complete intersection itself. 
The respective (Cox) ideals are
$$
  J_{\hat\Sigma} = \langle \prod_{i=0}^{3+\ell}x_i \rangle,
  \quad \mathrm{and}\qquad
  J_{\hat\Xi} = \langle \prod_{a=1}^{\ell}p_a \rangle.
$$
For the following let $\M = \M^1_{\hat\Sigma}$ and $\L := \M^1_{\hat\Xi}$ abbreviate the autoequivalences 
of twisting by $\hat S(1)$ in the two associated singularity categories.

\

\noindent\emph{In the category $\calc_{\hat\Sigma}$}:

\

\noindent For $(G_1,\ldots,G_\ell)$ generic, according to Proposition~\ref{twistmonoRelation}, a 
regular sequence $\cals_\cali$ of four elements, say $\underline{x}=(x_{i_0},x_{i_1},x_{i_2},x_{i_3})$, 
gives rise to a complex isomorphic to the zero object,
$$
  A[4] 
  ~~\stackrel{\underline{x}}{\longrightarrow}
  ~~\bigoplus_{b=0}^3\M^{w_{i_b}}(A)[3]
  ~~\stackrel{\underline{x}}{\longrightarrow}
  ~~\ldots
  ~~\stackrel{\underline{x}}{\longrightarrow}
  ~~\bigoplus_{b=0}^3\M^{w-w_{i_b}}(A)[1]
  ~~\stackrel{\underline{x}}{\longrightarrow}
  ~~\M^{w}(A) ~~\cong~~ 0,
$$
where $w=\sum_{b=0}^3 w_{i_b}$. 
In K-theory for each regular sequence $\cals_\cali$ the relation becomes $\prod_{b=0}^3(M^{w_{i_b}}-\id)=0$. 

As an example, for $X(10) \subset \bbP^4(11125)$ the K-theory relations are
\begin{eqnarray}
  \nonumber
  (M-\id)^2(M^2-\id)(M^5-\id) &=& 0,\\
  \nonumber
  (M-\id)^3(M^2-\id) &=& 0,\\
  \nonumber
  (M-\id)^3(M^5-\id) &=& 0.
\end{eqnarray}
It is easy to verify that for each complete intersection in the above list, the relations imply the 
well-known result that the autoequivalence $M$ is \emph{maximally} unipotent, that is
$$
  (M - \id)^4 = 0.
$$

\

\noindent\emph{In the category $\calc_{\hat\Xi}$}:

\

\noindent Proposition~\ref{twistmonoRelationPhi} tells us that for $\underline{p}=(p_1,\ldots,p_\ell)$ the complex
\begin{equation}
  \label{SVcomplex}
  \L^d(A^\cdot) 
  ~~\stackrel{\underline{p}}{\longrightarrow}
  ~~\oplus_a\L^{d-d_a}(A^\cdot)[1]
  ~~\stackrel{\underline{p}}{\longrightarrow}
  ~~\ldots
  ~~\stackrel{\underline{p}}{\longrightarrow}
  ~~\oplus_a\L^{d_a}(A^\cdot)[\ell\!-\!1]
  ~~\stackrel{\underline{p}}{\longrightarrow}
  ~~A^\cdot[\ell],
\end{equation}
is isomorphic to the zero object. Here, $d=\sum_a d_a$. 
On the level of K-theory this becomes
$$
  \prod_{a=1}^\ell(L^{d_a}-\id) = 0.
$$

\

\noindent\emph{Back in $\calc_{\hat\Sigma}$}:

\

\noindent Theorem~\ref{evenbigger} uses the functor $F_{\hat\Xi}:\calc_{\hat\Sigma}\longrightarrow 
\calc_{\hat\Xi}$ to map the relation (\ref{SVcomplex}) to a relation for the autoequivalence 
$\L_{F_{\hat\Xi}} = F_{\hat\Xi}^{-1} \L F_{\hat\Xi}$ on $\calc_{\hat\Sigma}$. 
Although the relation for $\L_{F_{\hat\Xi}}$ is a straight forward consequence of (\ref{SVcomplex}), 
we stress that it is highly non-trivial, even more so, if we ``forget'' that $\L_{F_{\hat\Xi}}$ in 
fact comes from an action of the Picard lattice in another phase. 

We illustrate this for the complete intersection $X=X(3,3)$ of two cubics in  $\bbP^5$, applying the 
autoequivalence $\L_{F_{\hat\Xi}}$ on the structure sheaf $\calo$ of $X$. Let $\Omega$ be the 
pull-back of the cotangent bundle of $\bbP^5$ to $X$. We obtain
$$
  \begin{array}{rccccccccc}
  \L_{F_{\hat\Xi}}(\calo)   &=& &&&&\Omega(1)[1]\ ,\\[10pt]
  (\L_{F_{\hat\Xi}}{})^2(\calo) &=& &&&&\wedge^2 \Omega(2)[2]\ ,\\[10pt]
  (\L_{F_{\hat\Xi}}{})^3(\calo) &=& &&\calo[2]^{\oplus 2}&\mapup{\varphi_3}&\wedge^3 \Omega(3)[3]\ ,\\[10pt]
  (\L_{F_{\hat\Xi}}{})^4(\calo) &=& &&\Omega(1)[3]^{\oplus 2}&\mapup{\varphi_4}&\wedge^4 \Omega(4)[4]\
,\\[10pt]
  (\L_{F_{\hat\Xi}}{})^5(\calo) &=& &&\wedge^2 \Omega(2)[4]^{\oplus 2}&\mapup{\varphi_5}&\wedge^5 \Omega(5)[5]\
,\\[10pt]
  (\L_{F_{\hat\Xi}}{})^6(\calo) &=&  \calo[4]^{\oplus 3}&\mapup{\varphi'_3}&\wedge^3 \Omega(3)[5]^{\oplus 2}
,
  \end{array}
$$\\
where the arrows are canonical elements in ${\rm Ext}^1(-,-)$.
In fact, the object in the last line is isomorphic to
$$
  \begin{array}{rcccc}
  (\L_{F_{\hat\Xi}}{})^6(\calo) &=& \calo[4]^{\oplus 4}&\mapup{\varphi_3}&\wedge^3 \Omega(3)[5]^{\oplus 2}, \\[-17pt]
  &&&\mapdiagdown{(1000)}&\oplus\\[5pt]
  &&&&\calo[3]
  \end{array}
$$
which confirms the relation for $\L_{F_{\hat\Xi}}$, when applied to $\calo$,
$$
    (\L_{F_{\hat\Xi}})^6(\calo) 
  ~~\cong~~
  ~~(\L_{F_{\hat\Xi}})^3(\calo)[2]^{\oplus 2}
  ~~\stackrel{F_{\hat\Xi}(\underline{p})}{\longrightarrow}
  ~~\calo[3].
$$

\end{document}